\documentclass{amsart}

\usepackage{amsmath}
\usepackage{amssymb}
\usepackage{amsthm}
\usepackage{mathabx}
\usepackage{amsrefs}        
\usepackage{pgfplots}
\usepackage{tikz}
\usetikzlibrary{trees}
\usetikzlibrary{arrows}
\usepackage{hyperref}
\usepackage[T2A]{fontenc}
\usepackage[utf8]{inputenc}
\usepackage{bbm}
\usepackage{bbold}
\usepackage{datetime}
\DeclareMathAlphabet{\mathbbold}{U}{BOONDOX-ds}{m}{n}
\usepackage{cleveref}

\newfont{\bbf}{msbm10 scaled\magstep1}

\def\arxivprefix{https://arxiv.org}
\newcommand*\arXiv[1]{\href{\arxivprefix/abs/#1}{arXiv:#1}}
\def\doiprefix{http://dx.doi.org/}

\DefineSimpleKey{bib}{primaryclass}{}
\DefineSimpleKey{bib}{archiveprefix}{}

\BibSpec{arXiv}{%
	+{}{\PrintAuthors}{author}
	+{,}{ \textit}{title}
	+{}{ \parenthesize}{date}
	+{,} { } {note}
	+{,}{ \arXiv }{eprint}
	+{,}{ primary class }{primaryclass}
}

\newcounter{glob}[section]
\renewcommand\theglob{%
	\ifnum\arabic{section}=0\else\arabic{section}.\fi %
	\arabic{glob}}

\theoremstyle{plain}

\crefname{thm}{theorem}{theorems}
\Crefname{thm}{Theorem}{Theorems}

\crefname{lemma}{lemma}{lemmas}
\Crefname{lemma}{Lemma}{Lemmas}

\crefname{cor}{corollary}{corollaries}
\Crefname{cor}{Corollary}{Corollaries}
\newtheorem{prop}[glob]{Proposition}
\crefname{prop}{proposition}{propositions}
\Crefname{prop}{Proposition}{Propositions}
\newtheorem*{prop*}{Proposition}

\theoremstyle{definition}
\newtheorem{defi}[glob]{Definition}
\crefname{defi}{definition}{definitions}
\Crefname{defi}{Definition}{Definitions}

\crefname{ex}{example}{examples}
\Crefname{ex}{Example}{Examples}

\crefname{question}{question}{questions}
\Crefname{question}{Question}{Questions}

\theoremstyle{remark}

\def\Z{\mathbb{Z}}
\def\ZZ{\mathbb{Z}}

\pgfplotsset{width=8cm}

\newdateformat{UKvardate}{%
	\THEDAY\ \monthname[\THEMONTH], \THEYEAR}
\UKvardate

\title{\bf Actions of frieze groups on inverse limits of polynomial rings}
\author{Elitza Hristova}
\address{Institute of Mathematics and Informatics,
	Bulgarian Academy of Sciences,
	Acad. G. Bonchev Str., Block 8,
	1113 Sofia, Bulgaria}
\email{e.hristova@math.bas.bg}
\author{Bogdan Stankov}
\address{Institute of Mathematics and Informatics,
	Bulgarian Academy of Sciences,
	Acad. G. Bonchev Str., Block 8,
	1113 Sofia, Bulgaria}
\email{bogdan.zl.stankov@gmail.com}
\date{\today}

\newcommand*{\thekeywords}{frieze groups, rings of invariants, symmetric functions}


\thanks{
The second author is supported by the Bulgarian Ministry of Education and Science, Scientific Programme “Enhancing the Research Capacity in Mathematical Sciences (PIKOM)”, No. DO1-241/15.08.2023.
}

\makeatletter
\hypersetup{
pdfauthor={Elitza Hristova and Bogdan Stankov},
pdftitle=\@title,
pdfsubject=\@title,
pdfkeywords=\thekeywords,
pdflang={en-US}
}
\makeatother 

\begin{document}

\maketitle

\begin{abstract}
In the spirit of the action of the symmetric group on the ring of polynomials in $n$ variables, we consider the actions of the seven frieze groups on rings of formal infinite linear combinations of monomials of restricted degree.
For each group we describe the respective subring of invariants.
We discuss also the structure of those rings as modules over each frieze group.
\end{abstract}

\section{Introduction}
\thispagestyle{empty}

Let $K$ be a field of characteristic $0$ an let $X = \{x_i\}_{i \in \ZZ}$ and $Y = \{y_i\}_{i \in \ZZ}$ be two infinite sets of variables indexed by the ring of integers $\ZZ$.
By $P_k[X]$ (respectively, $P_k[X,Y]$) we denote the ring of formal infinite linear combinations of monomials of degree $k$ from $X$ (respectively, from $X$ and $Y$) with coefficients in $K$.
We also take 
\[
P[X] = \bigoplus_{k \in \ZZ_{\geq 0}} P_k[X] \text{ and } P[X,Y] = \bigoplus_{k \in \ZZ_{\geq 0}} P_k[X,Y].
\] 
Then there is a classical action of the infinite symmetric group $S_{\infty}$ on the ring $P[X]$ and in particular on each homogeneous component $P_k[X]$ given by permuting the variables. The ring of symmetric functions is then the subring of invariants in $P[X]$ under the action of $S_{\infty}$. It is a classical problem to study the structure and properties of the ring of symmetric functions. In particular, several bases with many applications outside the theory of symmetric functions have been given like e.g., the elementary symmetric functions, the complete symmetric functions, the Schur functions. More information on this topic can be found in e.g., \cite{fultonharris}, \cite{macdonald98}, and many more.

In this paper, in analogy with the group $S_{\infty}$, we consider the action of the so-called frieze groups on the rings $P[X]$ and $P[X,Y]$.  

A two-dimensional figure or pattern in the plane which repeats infinitely in one direction is called a frieze pattern. The group of symmetries of a frieze pattern, i.e. the group of transformations in the plane which leave the frieze pattern invariant, is called a frieze group. In other words (see \cite{cederberg}), a frieze group is a group of transformations that keep a given line $c$ invariant and whose translations form an infinite cyclic subgroup.
Therefore, frieze groups are a class of infinite discrete groups. There are seven distinct frieze groups which we will denote by $F_i$ for $i = 1, \dots, 7$ (see, e.g., \cite{cederberg}, \cite{coxeter69} for the classification of frieze groups).

For the frieze groups $F_1$ and $F_3$ there is a natural action on the ring $P[X]$ and on each homogeneous component  $P_k[X]$. Similarly, for the other five frieze groups there is a natural action on the ring $P[X,Y]$ and respectively, on $P_k[X,Y]$.
The first goal of the current paper is to describe the subring of invariants of $P_k[X]$ and of $P_k[X,Y]$ under the action of each frieze group. In \Cref{sec_prelim} we give the necessary definitions of the seven frieze groups and their actions on $P_k[X]$ and on $P_k[X,Y]$.
Then in \Cref{sec_inv} we describe the respective subrings of invariants.
Finally, in \Cref{sec_Structure}, we use the results from the previous sections to describe the structure of $P_k[X]$ for each $k$ as a representation of $F_1$ and similarly the structure of $P_k[X,Y]$ for each $k$ as a representation of $F_6$. Similar results can be obtained for the other five frieze groups.

\section*{Acknowledgments}
The authors are grateful to Vesselin Drensky for suggesting this topic and for the stimulating discussions during the work on this text.

\section{Preliminaries}\label{sec_prelim}

\subsection{Frieze groups}
In this section we recall the definitions of the seven frieze groups.

The frieze group $F_1$ can be defined as the group of symmetries of the following sample pattern:

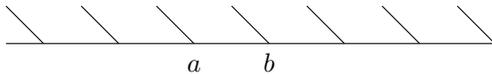
\begin{center}
	\vspace*{0.5cm}
	\begin{tikzpicture}
		\draw (0,0) -- (6.5,0);
		\draw (0, 0.5) -- (0.5, 0);
		\draw (1, 0.5) -- (1.5, 0);
		\draw (2, 0.5) -- (2.5, 0);
		\draw (2.5, -0.1) node[below] {$a$}; 
		\draw (3, 0.5) -- (3.5, 0);
		\draw (3.5, -0) node[below] {$b$};
		\draw (4, 0.5) -- (4.5, 0);
		\draw (5, 0.5) -- (5.5, 0);
		\draw (6, 0.5) -- (6.5, 0);
	\end{tikzpicture}
	\vspace*{0.5cm}
\end{center}

Similar patterns for $F_1$ and for the other frieze groups can be found on Wikipedia. Other sample patterns can be found, e.g., in \cite{cederberg} and \cite{coxeter69}.

$F_1$ is generated by a single translation by the distance between the points $a$ and $b$. Thus, it is isomorphic to the infinite cyclic group $\mathbb{Z}$. In terms of generators and relations, $F_1$ can be written as $F_1 = \left\langle t \right\rangle \cong \Z$.

The frieze group $F_2$ can be defined as the group of symmetries of the following sample pattern:

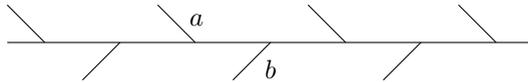
\begin{center}
	\vspace*{0.5cm}
	\begin{tikzpicture}
		\draw (0,0) -- (7,0);
		\draw (0, 0.5) -- (0.5, 0);
		\draw (1.5, 0) -- (1, -0.5);
		\draw (2, 0.5) -- (2.5, 0);
		\draw (2.3, 0.3) node[right]{$a$};
		\draw (3.5, 0) -- (3, -0.5);
		\draw (3.3, -0.35) node[right] {$b$};
		\draw (4, 0.5) -- (4.5, 0);
		\draw (5.5, 0) -- (5, -0.5);
		\draw (6, 0.5) -- (6.5, 0);
	\end{tikzpicture}
	\vspace*{0.5cm}
\end{center}

$F_2$ is generated by one glide reflection $g$ which sends the line segment $a$ to the line segment $b$. Therefore, $F_2$ is also isomorphic to $\ZZ$. In terms of generators and relations, $F_2$ can be written as
\[
F_2 = \left\langle g \right\rangle \cong \Z.
\]

It is useful to notice that $F_2$ contains also a translation $t'$, which is given by $t' = g^2$.

The frieze group $F_3$ can be defined as the group of symmetries of the following sample pattern:

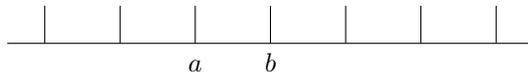
\begin{center}
	\vspace*{0.5cm}
	\begin{tikzpicture}
		\draw (0,0) -- (7,0);
		\draw (0.5, 0.5) -- (0.5, 0);
		\draw (1.5, 0.5) -- (1.5, 0);
		\draw (2.5, 0.5) -- (2.5, 0);
		\draw (2.5, -0.1) node[below] {$a$}; 
		\draw (3.5, 0.5) -- (3.5, 0);
		\draw (3.5, -0) node[below] {$b$};
		\draw (4.5, 0.5) -- (4.5, 0);
		\draw (5.5, 0.5) -- (5.5, 0);
		\draw (6.5, 0.5) -- (6.5, 0);
	\end{tikzpicture}
	\vspace*{0.5cm}
\end{center}

Therefore, $F_3$ is generated by one translation by the distance between the points $a$ and $b$ and one vertical reflection for example along the vertical line going through the point $a$. In terms of generators and relations, $F_3$ can be written as
\[
F_3 = \left\langle t, v : v^2 = 1, (vt)^2 = 1\right\rangle.
\]

Hence, $F_3$ is isomorphic to the infinite dihedral group $D_{\infty}$.

The frieze group $F_4$ can be defined as the group of symmetries of the following sample pattern:

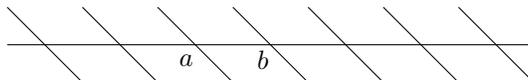
\begin{center}
	\vspace*{0.5cm}
	\begin{tikzpicture}
		\draw (0,0) -- (7,0);
		\draw (0, 0.5) -- (1, -0.5);
		\draw (1, 0.5) -- (2, -0.5);
		\draw (2, 0.5) -- (3, -0.5);
		\draw (2.6, -0.2) node[left] {$a$}; 
		\draw (3, 0.5) -- (4, -0.5);
		\draw (3.6, -0.2) node[left] {$b$};
		\draw (4, 0.5) -- (5, -0.5);
		\draw (5, 0.5) -- (6, -0.5);
		\draw (6, 0.5) -- (7, -0.5);
	\end{tikzpicture}
	\vspace*{0.5cm}
\end{center}
$F_4$ is generated by one translation by the distance between the points $a$ and $b$ and one $180^\circ$ rotation about for example the point $a$. In terms of generators and relations, $F_4$ can be written as
\[
F_4 = \left\langle t, r : r^2 = 1, (rt)^2 = 1 \right\rangle \cong D_{\infty}.
\]

The frieze group $F_5$ can be defined as the group of symmetries of the following sample pattern:

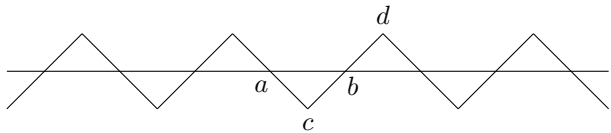
\begin{center}
	\vspace*{0.5cm}
	\begin{tikzpicture}
		\draw (-1,0) -- (7,0);
		\draw (-1, -0.5) -- (0, 0.5);
		\draw (0, 0.5) -- (1, -0.5);
		\draw (1, -0.5) -- (2, 0.5);
		\draw (2, 0.5) -- (3, -0.5);
		\draw (3, -0.5) -- (4, 0.5);
		\draw (2.6, -0.2) node[left] {$a$}; 
		\draw (3.4, -0.2) node[right] {$b$};
		\draw (3, -0.5) node[below] {$c$};
		\draw (4, 0.5) node[above] {$d$};
		\draw (4, 0.5) -- (5, -0.5);
		\draw (5, -0.5) -- (6, 0.5);
		\draw (6, 0.5) -- (7, -0.5);
	\end{tikzpicture}
	\vspace*{0.5cm}
\end{center}

There are four transformations which leave this pattern invariant. These are translation by the distance between the points $a$ and $b$; vertical reflection for example along the vertical line passing through the point $c$; $180^{\circ}$ rotation with center for example the point $a$ and a glide reflection which sends the point $c$ to the point $d$. Since some of these transformations can be expressed in terms of the others, using generators and relations the group $F_5$ can be defined in the following way:
\[
F_5 = \left\langle g, v: v^2 = 1, (vg)^2 = 1 \right\rangle \cong D_{\infty}.
\] 
Therefore, $F_5$ is generated by one glide reflection $g$ and one vertical reflection $v$. Then there is a translation $t'$, which can be expressed as $t' = g^2$ and a rotation $r'$, which can be expressed as $r' = v g$.

The frieze group $F_6$ can be defined as the group of symmetries of the following sample pattern:

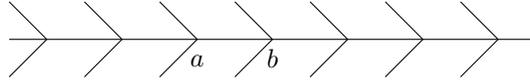
\begin{center}
	\vspace*{0.5cm}
	\begin{tikzpicture}
		\draw (0,0) -- (7,0);
		\draw (0, 0.5) -- (0.5, 0);
		\draw (0.5, 0) -- (0, -0.5);
		\draw (1, 0.5) -- (1.5, 0);
		\draw (1.5, 0) -- (1, -0.5);
		\draw (2, 0.5) -- (2.5, 0);
		\draw (2.5, 0) -- (2, -0.5);
		\draw (2.5, -0.08) node[below] {$a$}; 
		\draw (3, 0.5) -- (3.5, 0);
		\draw (3.5, 0) -- (3, -0.5);
		\draw (3.5, 0) node[below] {$b$};
		\draw (4, 0.5) -- (4.5, 0);
		\draw (4.5, 0) -- (4, -0.5);
		\draw (5, 0.5) -- (5.5, 0);
		\draw (5.5, 0) -- (5, -0.5);
		\draw (6, 0.5) -- (6.5, 0);
		\draw (6.5, 0) -- (6, -0.5);
	\end{tikzpicture}
	\vspace*{0.5cm}
\end{center}

$F_6$ is generated by one translation by the distance between the points $a$ and $b$ and one reflection along the horizontal axis. In terms of generators and relations, $F_6$ can be written as
\[
F_6 = \left\langle t, h: h^2 = 1, th=ht \right\rangle \cong \Z\times \Z/2\Z.
\] 

The frieze group $F_7$ can be defined as the group of symmetries of the following sample pattern:

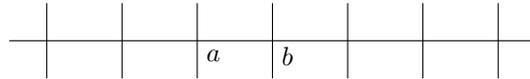
\begin{center}
	\vspace*{0.5cm}
	\begin{tikzpicture}
		\draw (0,0) -- (7,0);
		\draw (0.5, 0.5) -- (0.5, -0.5);
		\draw (1.5, 0.5) -- (1.5, -0.5);
		\draw (2.5, 0.5) -- (2.5, -0.5);
		\draw (2.5, -0.2) node[right] {$a$}; 
		\draw (3.5, 0.5) -- (3.5, -0.5);
		\draw (3.5, -0.2) node[right] {$b$};
		\draw (4.5, 0.5) -- (4.5, -0.5);
		\draw (5.5, 0.5) -- (5.5, -0.5);
		\draw (6.5, 0.5) -- (6.5, -0.5);
	\end{tikzpicture}
	\vspace*{0.5cm}
\end{center}

The group $F_7$ requires three generators. One generating set consists of a translation $t$ by the distance between the points $a$ and $b$, a reflection $h$ across the horizontal axis and one vertical reflection $v$ for example across the vertical line going through the point $a$. Then, there is one more transformation which leaves this pattern invariant, namely a $180^{\circ}$ rotation $r$ with center for example the point $a$. This rotation $r$ can be expressed as $r = hv = vh$. Therefore, in terms of generators and relations, $F_7$ can be written as
\[
F_7 = \left\langle t, v, h : v^2 = 1, (vt)^2 = 1, h^2 = 1, th = ht, vh = hv \right\rangle.
\]

Hence, $F_7$ is isomorphic to the group $D_{\infty} \times \Z/2\Z$.

\subsection{Actions on the rings $P[X]$ and $P[X,Y]$}
Formally, the rings $P[X]$ and $P[X,Y]$ are defined in the following way.
Let $K[x_{-n}, \dots, x_{n}]$ be the ring of polynomials in $2n+1$ variables with coefficients from $K$ and let $K_k[x_{-n}, \dots, x_{n}]$ denote its homogeneous component of total degree $k$. For any $m \geq n$ there is a natural homomorphism
\[
\rho_{m,n}^k: K_k[x_{-m}, \dots, x_m] \rightarrow K_k[x_{-n}, \dots, x_n],
\]
which sends the variables $x_{n+1}, \dots, x_m$ and $x_{-m}, \dots, x_{-n-1}$ to zero and the variables $x_{-n}, \dots, x_n$ to themselves. We take the inverse limit
\[
P_k[X] = \varprojlim_n K_k[x_{-n}, \dots, x_n]
\] 
relative to the homomorphisms $\rho_{m,n}^k$ and set 
\[
P[X]= \bigoplus_{k\geq 0} P_k[X].
\]
Thus, every element of $P_k[X]$ is a formal infinite linear combination of monomials $x_{i_1}\cdots x_{i_k}$. 

The rings $P_k[X,Y]$ for $k \geq 1$ are defined in a similar way using two sets of variables:
\[
P_k[X,Y] = \varprojlim_n K_k[x_{-n}, \dots, x_n; y_{-n} \dots y_{n}]
\] 
relative to the homomorphisms $\eta_{m,n}^k$ which are defined as
\[
\eta_{m,n}^k: K_k[x_{-m}, \dots, x_m; y_{-m}, \dots, y_m] \rightarrow K_k[x_{-n}, \dots, x_n; y_{-n} \dots y_n],
\]
which sends the variables $x_{n+1}, \dots, x_m; y_{n+1}, \dots, y_m$ and $x_{-m}, \dots, x_{-n-1}; y_{-m} \dots y_{-n-1}$ to zero and the other variables to themselves. Then the ring $P[X,Y]$ is defined as
\[
P[X,Y] = \bigoplus_{k\geq 0} P_k[X,Y].
\]

The group $F_1$ acts on the ring of integers $\mathbb{Z}$ in the following natural way $t . z = z+1$ for any $z \in \ZZ$. This induces an action of $F_1$ on the ring $P_k[X]$ for any $k \geq 1$. The action is defined on the monomials in the following way:
\begin{align*}
	& t: x_i \mapsto x_{t.i} = x_{i+1}\\
	& t: x_{i_1}\cdots x_{i_k} \mapsto x_{t.i_1} \cdots x_{t. i_k} = x_{i_1 +1 } \cdots x_{i_k +1}.
\end{align*}
 
Then, the action of $F_1$ defined above on the set of monomials is extended by linearity to any element of $P_k[X]$.

Similarly, we can define actions of $F_3$ on $P_k[X]$ and of the other five frieze groups on $P_k[X,Y]$. Below we give the respective formulas.

The group $F_2$ acts naturally on the ring $P_k[X,Y]$ for any $k \geq 1$ and this action can be described on the set of monomials in the following way:
\begin{align*}
	& g: x_{i} \mapsto y_{i+1}\\
	& g: y_i \mapsto x_{i+1}.
\end{align*} 
Then, we obtain that the translation $t'$ acts on the monomials by the formula $t'(x_i) = g^2(x_i) = x_{i+2}$ and similarly $t'(y_i) = y_{i+2}$.
Remark that the action of $t'$ on  functions containing only variables in $\{x_i\}_{i \in \ZZ}$ is the same as the action of $t^2$ for the translation $t$ defined in $F_1$.
It is more generally the same as the action of $t^2$ for the translation $t$ in $F_3$, $F_4$, $F_6$ and $F_7$.
Notice however that the group $F_2$ does not contain the translation $t$ - only the translation of double length $t'$.

The action of $F_3$ on the ring $P_k[X]$ for any $k \geq 1$ can again be given on the monomials in the following way:
\begin{align*}
	& t: x_{i} \mapsto x_{i + 1} \\
	& v: x_{i} \mapsto x_{-i}.
\end{align*} 

The group $F_4$ acts on the ring $P_k[X,Y]$ for any $k \geq 1$ and we again describe this action on the monomials:
\begin{align*}
	& t: x_i \mapsto x_{i+1}\\
	& t: y_i \mapsto y_{i+1}\\
	&\\
	& r: x_i \mapsto y_{-i}\\
	& r : y_i \mapsto x_{-i}.
\end{align*}

The group $F_5$ acts again on $P_k[X,Y]$ for each $k \geq 0$ and this action is given by:
\begin{align*}
	& g: x_i \mapsto y_{i+1}\\
	& g: y_i \mapsto x_{i+1}\\
	&\\
	& v: x_i \mapsto x_{-i}\\
	& v: y_i \mapsto y_{-i}.
\end{align*}
Then, for the translation $t'$, as in the case of $F_2$, we obtain $t'(x_i) = g^2(x_i) = x_{i+2}$ and $t'(y_i) = y_{i+2}$. In the same way, for the rotation $r'$ we obtain $r'(x_i) = vg(x_i) = y_{-i-1}$ and similarly $r'(y_i) = x_{-i-1}$. 	

The action of $F_6$ on $P_k[X,Y]$ is given by:
\begin{align*}	
	& t: x_i \mapsto x_{i+1}\\
	& t: y_i \mapsto y_{i+1}\\
	&\\
	& h: x_i \mapsto y_i\\
	& h: y_i \mapsto x_i.
\end{align*}

Finally, the action of $F_7$ on $P_k[X,Y]$ is given by:
\begin{align*}
	& t: x_i \mapsto x_{i+1}\\
	& t: y_i \mapsto y_{i+1}\\
	&\\
	& v: x_i \mapsto x_{-i}\\
	& v: y_i \mapsto y_{-i}\\
	&\\
	& h: x_i \mapsto y_i\\
	& h: y_i \mapsto x_i.
\end{align*}

Then, the action of the rotation $r$ is obtained by $r(x_i) = hv(x_i) = y_{-i}$ and $r(y_i) = x_{-i}$.

\section{Rings of invariants} \label{sec_inv}

\subsection{The frieze group $F_1$} \label{sec_F1}

Our first goal is to describe the invariant elements of $P_k[X]$ under the action of $F_1$.

Notice that if a function of the form $\sum_{i_1\leq i_2\leq\dots\leq i_k\in\Z^k}c_{i_1,i_2,\dots,i_k}x_{i_1}x_{i_2}\dots x_{i_k}$ is invariant for the action of a group $G$, then for every $g\in G$, $c_{i_1,i_2,\dots,i_k}=c_{g.i_1,g.i_2,\dots,g.i_k}$.
It is not hard to check that this condition is also sufficient.
Therefore we need to understand the orbits of the induced action of $F_1$ on weakly increasing sequences in $\Z^k$ for any given fixed $k$.

Consider $i_1\leq i_2\leq\dots\leq i_k$ and $g\in F_1$.
We then have $t^{1-i_1}(i_1,\dots,i_k)=(1,i_2-i_1 +1,\dots,i_k-i_1 +1)$.
For any $g\in F_1$, $g\neq t^{1-i_1}$, we have $g=t^z$ for some $z\in\Z$ (and thus $z\neq1-i_1$).
Then $z(i_1,\dots,i_k)=(i_1+z,\dots,i_k+z)$ with $i_1+z\neq1$.
Therefore in each orbit, there is exactly one sequence such that the smallest number is $1$.
We describe such a monomial using the sequence of its powers:

\begin{defi}\label{composition}
A \textit{non-negative composition} of order $k\in\Z_{>0}$ we will call a sequence $\sigma=(\sigma_1,\dots,\sigma_m)$ of non-negative integers such that $\sigma_1>0$, $\sigma_m>0$ and

$$\sum_{i=1}^m\sigma_i=k.$$

We denote the order as $ord(\sigma)=k$ and the number of parts as $par(\sigma)=m$.

Additionally, we will also consider the empty sequence $\sigma=\emptyset$ as a non-negative composition with $ord(\emptyset)=0$ and $par(\emptyset)=0$.
\end{defi}

Remark that a non-negative composition of fixed order $k$ for $k\geq2$ can contain any number of zeros, and therefore there are infinitely many non-negative compositions of that order.
Whereas for $k=1$, we necessarily have $m=1$ and $\sigma_1=\sigma_m=1>0$.

We then have that for any monomial $X=x_{i_1}x_{i_2}\dots x_{i_k}$ such that the smallest number $j$ for which $x_j$ is present is $j=1$ can be written in a unique way as 

$$X=x_1^{\sigma_1}x_2^{\sigma_2}\dots x_{par(\sigma)}^{\sigma_{par(\sigma)}}$$
where $\sigma$ is a non-negative composition.
Similarly, each monomial in $P_k[X]$ can be written in a unique way as $x_{i+1}^{\sigma_1}x_{i+2}^{\sigma_2}\dots x_{i+par(\sigma)}^{\sigma_{par(\sigma)}}$ for some $i\in\Z$.

If $\sigma$ is the empty non-negative composition, this monomial will be evaluated as $1$.
Using that convention, a monomial in $P_k[X,Y]$ can then written in a unique way in the form $x_{i+1}^{\sigma_1}x_{i+2}^{\sigma_2}\dots x_{par(\sigma)}^{\sigma_{par(\sigma)}}y_{i+1+\Delta}^{\sigma'_1}y_{i+2+\Delta}^{\sigma'_2}\dots y_{i+par(\sigma')+\Delta}^{\sigma'_{par(\sigma')}}$ with $ord(\sigma)+ord(\sigma')=k$ (where one of those orders might be zero) and $i,\Delta \in \ZZ$ (if $\sigma$ or $\sigma'$ is of order $0$, $\Delta$ is taken to be $0$).

We obtain that orbits of $F_1$ are of the form $t^z.x_1^{\sigma_1}x_2^{\sigma_2}\dots x_{par(\sigma)}^{\sigma_{par(\sigma)}}$ where $\sigma$ is fixed and $z$ varies in $\Z$.
By summing over an orbit, we obtain invariant functions:

\begin{equation*}\label{f1funct}
f^{(1)}_{\sigma}=\sum_{i\in\Z}t^i.x_1^{\sigma_1}x_2^{\sigma_2}\dots x_{par(\sigma)}^{\sigma_{par(\sigma)}}.
\end{equation*}

As each monomial in $P_k[X]$ is present in exactly one of those functions, we obtain:

\begin{prop}
For $k\geq1$, each element of $P_k[X]^{F_1}$ can be written in a unique way as a formal infinite linear combination on the set
$$\left\{f^{(1)}_\sigma|\sigma\mbox{ is a non-negative composition with }ord(\sigma)=k\right\},$$
where
$$f^{(1)}_\sigma=\sum_{i\in\Z}t^i.x_1^{\sigma_1}x_2^{\sigma_2}\dots x_{par(\sigma)}^{\sigma_{par(\sigma)}}=\sum_{i\in\Z}x_{i+1}^{\sigma_1}x_{i+2}^{\sigma_2}\dots x_{i+par(\sigma)}^{\sigma_{par(\sigma)}}.$$
\end{prop}

Let us compare this space to the space of symmetric functions (see~\cite[Chapter~II.2]{macdonald98}).
The translation $t$ has infinite support and is therefore not contained in the infinite symmetric group.
However, on any finite monomial, there exists an element of $S_n$ for large enough $n$ that acts on that monomial the same way as $t$.
It follows that any symmetric function is also invariant by the action of $F_1$; and therefore can be written as a formal infinite linear combination on the functions $f^{(1)}_\sigma$.
We will present this explicitly for the elementary symmetric functions $e_r$ and for the complete symmetric functions $h_r$ which give two bases for the space of all symmetric functions.
By definition (\cite{macdonald98}), the $r$-th elementary symmetric function $e_r$ for any $r \geq 0$ is defined as the formal infinite sum of all products of $r$ distinct variables $x_i$
\[
e_r = \sum_{i_1 < i_2 < \dots <i_r} x_{i_1}x_{i_2}\cdots x_{i_r}.
\] 
Thus, $e_r \in P_r[X]$, i.e., it is a homogeneous function of total degree $r$. It is not difficult to see from the definition that 

\begin{align}\label{eq_elem_func}
	e_r = \sum_{\sigma:ord(\sigma)=r,\sigma_i\leq1\forall i} f^{(1)}_\sigma.
\end{align}

Similarly, by definition (\cite{macdonald98}), the $r$-th complete symmetric function $h_r$ is the sum of all monomials of total degree $r$ in the variables $x_i$, $i \in \ZZ$, i.e.,
\[
h_r = \sum_{i_1 \leq i_2 \leq \dots  \leq i_r} x_{i_1}x_{i_2}\cdots x_{i_r}.
\] 

We have again that $h_r \in P_r[X]$. Furthermore,
\begin{align} \label{eq_complete_func}
	h_r = \sum_{\sigma:ord(\sigma)=r} f^{(1)}_\sigma.
\end{align}

In particular, \Cref{eq_elem_func,eq_complete_func} imply that for the Schur functions corresponding to the partition $(r, 0, \dots, 0)$ and $(\underbrace{1,\dots, 1}_r, 0, \dots, 0)$ we also have explicit expressions in terms of the functions $f^{(1)}_\sigma$.
Namely,

 \begin{align*}
 s_{(r, 0, \dots, 0)}=	h_r = \sum_{\sigma:ord(\sigma)=r} f^{(1)}_\sigma.
 \end{align*}

\begin{align*}
s_{(\underbrace{1,\dots, 1}_r, 0, \dots, 0)} = 	e_r = \sum_{\sigma:ord(\sigma)=r,\sigma_i\leq1\forall i} f^{(1)}_\sigma.
\end{align*}

\subsection{The frieze group $F_2$}\label{sec_F2}

Recall that this group contains a translation $t'=g^2$.
Furthermore, this translation acts the same way as $t^2$ for the transformation $t$ that appears $F_1$.
Therefore from the calculations in the previous section it follows that
$${t'}^i.x_1^{\sigma_1}x_2^{\sigma_2}\dots x_{par(\sigma)}^{\sigma_{par(\sigma)}}=x_{2i+1}^{\sigma_1}x_{2i+2}^{\sigma_2}\dots x_{2i+par(\sigma)}^{\sigma_{par(\sigma)}}$$
and similarly ${t'}^i.y_1^{\sigma_1}y_2^{\sigma_2}\dots y_{par(\sigma)}^{\sigma_{par(\sigma)}}=y_{2i+1}^{\sigma_1}y_{2i+2}^{\sigma_2}\dots y_{2i+par(\sigma)}^{\sigma_{par(\sigma)}}$.

Recall that a monomial in $P_k[X,Y]$ can be written in a unique way in the form $x_{i+1}^{\sigma_1}x_{i+2}^{\sigma_2}\dots x_{par(\sigma)}^{\sigma_{par(\sigma)}}y_{i+1+\Delta}^{\sigma'_1}y_{i+2+\Delta}^{\sigma'_2}\dots y_{i+par(\sigma')+\Delta}^{\sigma'_{par(\sigma')}}$ with $ord(\sigma)+ord(\sigma')=k$ and $i,\Delta \in \ZZ$ (if $\sigma$ or $\sigma'$ is of order $0$, $\Delta$ is taken to be $0$).
Using the previous calculations we get (we note $m=par(\sigma)$ and $m'=par(\sigma')$):

\begin{equation}\label{tp_twovar_eq}
{t'}^i.x_1^{\sigma_1}\dots x_{m}^{\sigma_{m}}y_{1+\Delta}^{\sigma'_1}\dots y_{m'+\Delta}^{\sigma'_{m'}}=x_{2i+1}^{\sigma_1}\dots x_{2i+m}^{\sigma_{m}}y_{2i+1+\Delta}^{\sigma'_1}\dots y_{2i+m'+\Delta}^{\sigma'_{m'}}
\end{equation}
and then

\begin{align*}
&g.x_{i+1}^{\sigma_1}x_{i+2}^{\sigma_2}\dots x_{i+m}^{\sigma_{m}}y_{i+1+\Delta}^{\sigma'_1}y_{i+2+\Delta}^{\sigma'_2}\dots y_{i+m'+\Delta}^{\sigma'_{m'}}\\
&=y_{i+1+1}^{\sigma_1}y_{i+1+2}^{\sigma_2}\dots y_{i+1+m}^{\sigma_{m}}x_{i+1+1+\Delta}^{\sigma'_1}x_{i+1+2+\Delta}^{\sigma'_2}\dots x_{i+1+m'+\Delta}^{\sigma'_{m'}}\\
&=x_{j+1}^{\sigma'_1}x_{j+2}^{\sigma'_2}\dots x_{j+m'}^{\sigma'_{m'}}y_{j+1-\Delta}^{\sigma_1}y_{j+2-\Delta}^{\sigma_2}\dots y_{j+m-\Delta}^{\sigma_{m}}
\end{align*}
for $j=i+1+\Delta$.

We consider if it is possible for a monomial to be fixed by $g^z$ for some $z\in\Z$.
As $t'=g^2$ and its powers are translations, they do not fix any monomials.
Therefore we have to consider and odd $z$.
From the calculation it follows that we need to have $\sigma=\sigma'$ and $\Delta=0$.
However, if $\Delta=0$, then the smallest index $j$ such that $x_j$ is present in the monomial changes parity when multiplied by $g^z$ for odd $z$.
Therefore

\begin{prop}\label{prop_f2}
For $k\geq1$, each element of $P_k[X]^{F_2}$ can be written in a unique way as a formal infinite linear combination on the set (see \Cref{composition})
$$\left\{f^{(2)}_{\sigma,\sigma',\Delta}|ord(\sigma)+ord(\sigma')=k,\Delta\in\Z\right\}\bigcup\left\{f^{(2)'}_{\sigma,\sigma',\Delta}|ord(\sigma)+ord(\sigma')=k,\Delta\in\Z\right\}$$
where

\begin{align*}
f^{(2)}_{\sigma,\sigma',\Delta} & =\sum_{i\in\Z}g^i.x_1^{\sigma_1}\dots x_{par(\sigma)}^{\sigma_{par(\sigma)}}y_{1+\Delta}^{\sigma'_1}\dots y_{par(\sigma')+\Delta}^{\sigma'_{par(\sigma')}} \\
& =\sum_{i\in\Z}x_{2i+1}^{\sigma_1}x_{2i+2}^{\sigma_2}\dots x_{2i+par(\sigma)}^{\sigma_{par(\sigma)}}y_{2i+1+\Delta}^{\sigma'_1}y_{2i+2+\Delta}^{\sigma'_2}\dots y_{2i+par(\sigma')+\Delta}^{\sigma'_{par(\sigma')}}\\
& +\sum_{i\in\Z}x_{2i+2+\Delta}^{\sigma'_1}x_{2i+3+\Delta}^{\sigma'_2}\dots x_{2i+par(\sigma')+1+\Delta}^{\sigma'_{par(\sigma')}}y_{2i+2}^{\sigma_1}y_{2i+3}^{\sigma_2}\dots y_{2i+par(\sigma)+1}^{\sigma_{par(\sigma)}}
&\\
f^{(2)'}_{\sigma,\sigma',\Delta} & =\sum_{i\in\Z}g^i.x_0^{\sigma_1}\dots x_{par(\sigma)-1}^{\sigma_{par(\sigma)}}y_{\Delta}^{\sigma'_1}\dots y_{par(\sigma')-1+\Delta}^{\sigma'_{par(\sigma')}} \\
& =\sum_{i\in\Z}x_{2i}^{\sigma_1}x_{2i+1}^{\sigma_2}\dots x_{2i+par(\sigma)-1}^{\sigma_{par(\sigma)}}y_{2i+\Delta}^{\sigma'_1}y_{2i+1+\Delta}^{\sigma'_2}\dots y_{2i+par(\sigma')-1+\Delta}^{\sigma'_{par(\sigma')}}\\
& +\sum_{i\in\Z}x_{2i+1+\Delta}^{\sigma'_1}x_{2i+2+\Delta}^{\sigma'_2}\dots x_{2i+par(\sigma')+\Delta}^{\sigma'_{par(\sigma')}}y_{2i+1}^{\sigma_1}y_{2i+2}^{\sigma_2}\dots y_{2i+par(\sigma)}^{\sigma_{par(\sigma)}}.
\end{align*}

Furthermore, for $\sigma=\sigma'$ and $\Delta=0$ we get $f^{(2)}_{\sigma,\sigma,0}=f^{(2)'}_{\sigma,\sigma,0}$.
Additionally, if $ord(\sigma)=0$ or $ord(\sigma')=0$, then $f^{(2)}_{\sigma,\sigma',\Delta}=f^{(2)}_{\sigma,\sigma',0}$ and $f^{(2)'}_{\sigma,\sigma',\Delta}=f^{(2)'}_{\sigma,\sigma',0}$ for all $\Delta\in\Z$.
\end{prop}

\subsection{The frieze group $F_3$} \label{sec_F3}

Remark that any element of $F_3$ can be written in the form $vt^z$ or $t^z$ for some $z\in\Z$.
Indeed $(vt)^2=1$ and therefore $tv=v^{-1}t^{-1}=vt^{-1}$ and $t^{-1}v=(vt)^{-1}=vt$, so we can always move $v$'s to the left in a word on those two generators.
And since $v^2=1$ we are left with $0$ or $1$ consecutive $v$'s.

The action of $t$ on monomials has been calculated in \Cref{sec_F1}, so we calculate the action of $v$:

\begin{equation} \label{calc_v_onevar}
v.x_{i+1}^{\sigma_1}\dots x_{i+par(\sigma)}^{\sigma_{par(\sigma)}} = x_{-i-1}^{\sigma_1}x_{-i-2}^{\sigma_2}\dots x_{-i-par(\sigma)}^{\sigma_{par(\sigma)}} = x_{j+1}^{\sigma_{par(\sigma)}}x_{j+2}^{\sigma_{par(\sigma)-1}}\dots x_{j+par(\sigma)}^{\sigma_1}
\end{equation}
where $j=-i-par(\sigma)-1$.

The orbit of an element $x_1^{\sigma_1}x_2^{\sigma_2}\dots x_{par(\sigma)}^{\sigma_{par(\sigma)}}$ by $F_3$ then consists of elements of the form $x_{i+1}^{\sigma_1}x_{i+2}^{\sigma_2}\dots x_{i+par(\sigma)}^{\sigma_{par(\sigma)}}$ or $x_{i+1}^{\sigma_{par(\sigma)}}x_{i+2}^{\sigma_{par(\sigma)-1}}\dots x_{i+par(\sigma)}^{\sigma_1}$.
Clearly, two elements of that form are distinct if defined by different $i$.
Therefore the only case in which they are not pairwise distinct is if $x_{i+1}^{\sigma_1}\dots x_{i+par(\sigma)}^{\sigma_{par(\sigma)}}=x_{i+1}^{\sigma_{par(\sigma)}}\dots x_{i+par(\sigma)}^{\sigma_1}$, which is the case if and only if $(\sigma_1,\sigma_2,\dots,\sigma_{par(\sigma)})=(\sigma_{par(\sigma)},\sigma_{par(\sigma)-1},\dots,\sigma_1)$.

\begin{prop}
For $k\geq1$, each element of $P_k[X]^{F_3}$ can be written in a unique way as a formal infinite linear combination on the set (see \Cref{composition})
$$\left\{f^{(3)}_\sigma|\sigma\mbox{ is a non-negative composition with }ord(\sigma)=k\right\},$$
where for $(\sigma_1,\dots,\sigma_{par(\sigma)})\neq(\sigma_{par(\sigma)},\dots,\sigma_1)$:

\begin{equation*}
\begin{split}
f^{(3)}_\sigma & = \sum_{i\in\Z}t^i.x_1^{\sigma_1}\dots x_{par(\sigma)}^{\sigma_{par(\sigma)}}+v.\sum_{i\in\Z}t^i.x_1^{\sigma_1}\dots x_{par(\sigma)}^{\sigma_{par(\sigma)}} \\
	& = \sum_{i\in\Z}x_{i+1}^{\sigma_1}x_{i+2}^{\sigma_2}\dots x_{i+par(\sigma)}^{\sigma_{par(\sigma)}}+\sum_{i\in\Z}x_{i+1}^{\sigma_{par(\sigma)}}x_{i+2}^{\sigma_{par(\sigma)-1}}\dots x_{i+par(\sigma)}^{\sigma_1}
\end{split}
\end{equation*}
and for $(\sigma_1,\dots,\sigma_{par(\sigma)})=(\sigma_{par(\sigma)},\dots,\sigma_1)$:

$$f^{(3)}_\sigma=\sum_{i\in\Z}t^i.x_1^{\sigma_1}x_2^{\sigma_2}\dots x_{par(\sigma)}^{\sigma_{par(\sigma)}}=\sum_{i\in\Z}x_{i+1}^{\sigma_1}x_{i+2}^{\sigma_2}\dots x_{i+par(\sigma)}^{\sigma_{par(\sigma)}}.$$

Furthermore, $f^{(3)}_\sigma=f^{(3)}_{(\sigma_{par(\sigma)},\dots,\sigma_1)}$ for all $\sigma$.

\end{prop}

\subsection{The frieze group $F_4$} \label{sec_F4}

Recall that $F_4\cong D_{\infty}\cong F_3$.
Equivalently to the beginning of \Cref{sec_F3}, we obtain that any element of $F_4$ can be written in the form $rt^z$ or $t^z$ for some $z\in\Z$.

We recall that a monomial in $P_k[X,Y]$ can be written in a unique way in the form $x_{i+1}^{\sigma_1}x_{i+2}^{\sigma_2}\dots x_{par(\sigma)}^{\sigma_{par(\sigma)}}y_{i+1+\Delta}^{\sigma'_1}y_{i+2+\Delta}^{\sigma'_2}\dots y_{i+par(\sigma')+\Delta}^{\sigma'_{par(\sigma')}}$ with $ord(\sigma)+ord(\sigma')=k$ and $i,\Delta \in \ZZ$ (if $\sigma$ or $\sigma'$ is of order $0$, $\Delta$ is taken to be $0$).
Similarly to \Cref{tp_twovar_eq}, we calculate (with $m=par(\sigma)$ and $m'=par(\sigma')$):

\begin{equation}\label{t_twovar_eq}
t^i.x_1^{\sigma_1}\dots x_{m}^{\sigma_{m}}y_{1+\Delta}^{\sigma'_1}\dots y_{m'+\Delta}^{\sigma'_{m'}}=x_{i+1}^{\sigma_1}\dots x_{i+m}^{\sigma_{m}}y_{i+1+\Delta}^{\sigma'_1}\dots y_{i+m'+\Delta}^{\sigma'_{m'}}.
\end{equation}

For $r$ we calculate:

\begin{equation}\label{rcalc}
\begin{split}
& r.x_{i+1}^{\sigma_1}x_{i+2}^{\sigma_2}\dots x_{i+m}^{\sigma_{m}}y_{i+1+\Delta}^{\sigma'_1}y_{i+2+\Delta}^{\sigma'_2}\dots y_{i+m'+\Delta}^{\sigma'_{m'}} \\
& = y_{j+1}^{\sigma_{m}}y_{j+2}^{\sigma_{m-1}}\dots y_{j+m}^{\sigma_1}x_{j+1-\Delta'}^{\sigma'_{m'}}x_{j+2-\Delta'}^{\sigma'_{m'-1}}\dots x_{j+m'-\Delta'}^{\sigma'_1} \\
& = x_{j'+1}^{\sigma'_{m'}}x_{j'+2}^{\sigma'_{m'-1}}\dots x_{j'+m'}^{\sigma'_1}y_{j'+1+\Delta'}^{\sigma_{m}}y_{j'+2+\Delta'}^{\sigma_{m-1}}\dots y_{j'+m+\Delta'}^{\sigma_1}
\end{split}
\end{equation}
for $j=-i-m-1$, $\Delta'=\Delta+m'-m$ and $j'=j-\Delta'$.

In order for those monomials to be pairwise distinct, we need to consider the case when $rt^z$ fixes a monomial $x_{i+1}^{\sigma_1}x_{i+2}^{\sigma_2}\dots x_{i+par(\sigma)}^{\sigma_{par(\sigma)}}y_{i+1+\Delta}^{\sigma'_1}y_{i+2+\Delta}^{\sigma'_2}\dots y_{i+par(\sigma')+\Delta}^{\sigma'_{par(\sigma')}}$ for some $z$.
From \Cref{rcalc}, that happens if and only if $(\sigma_1,\dots,\sigma_{par(\sigma)})=(\sigma'_{par(\sigma')},\dots,\sigma'_1)$ (remark that in that case $par(\sigma)=par(\sigma')$ and thus $\Delta'=\Delta$).

\begin{prop}
For $k\geq1$, each element of $(P_k[X,Y])^{F_4}$ can be written in a unique way as a formal infinite linear combination on the set (see \Cref{composition})
$$\left\{f^{(4)}_{\sigma,\sigma',\Delta}|\sigma,\sigma'\mbox{ are non-negative compositions with }ord(\sigma)+ord(\sigma')=k,\Delta\in\Z\right\}$$
where for $(\sigma_1,\dots,\sigma_{par(\sigma)})\neq(\sigma'_{par(\sigma')},\dots,\sigma'_1)$:

\begin{equation*}
\begin{split}
f^{(4)}_{\sigma,\sigma',\Delta} & = \sum_{i\in\Z}t^i.x_{1}^{\sigma_1}\dots x_{m}^{\sigma_{m}}y_{1+\Delta}^{\sigma'_1}\dots y_{m'+\Delta}^{\sigma'_{m'}}+r.\sum_{i\in\Z}t^i.x_{1}^{\sigma_1}\dots x_{m}^{\sigma_{m}}y_{1+\Delta}^{\sigma'_1}\dots y_{m'+\Delta}^{\sigma'_{m'}} \\
& = \sum_{i\in\Z}x_{i+1}^{\sigma_1}x_{i+2}^{\sigma_2}\dots x_{i+m}^{\sigma_{m}}y_{i+1+\Delta}^{\sigma'_1}y_{i+2+\Delta}^{\sigma'_2}\dots y_{i+m'+\Delta}^{\sigma'_{m'}}\\
&+\sum_{i\in\Z}x_{i+1}^{\sigma'_{m'}}x_{i+2}^{\sigma'_{m'-1}}\dots x_{i+m'}^{\sigma'_1}y_{i+1+\Delta'}^{\sigma_{m}}y_{i+2+\Delta'}^{\sigma_{m-1}}\dots y_{i+m+\Delta'}^{\sigma_1}
\end{split}
\end{equation*}
with $\Delta'=\Delta+par(\sigma')-par(\sigma)$, $m=par(\sigma)$ and $m'=par(\sigma')$; and for $(\sigma_1,\dots,\sigma_{par(\sigma)})=(\sigma'_{par(\sigma')},\dots,\sigma'_1)$:

\begin{align*}
f^{(4)}_{\sigma,\sigma',\Delta} & = \sum_{i\in\Z}t^i.x_{1}^{\sigma_1}\dots x_{par(\sigma)}^{\sigma_{par(\sigma)}}y_{1+\Delta}^{\sigma'_1}\dots y_{par(\sigma')+\Delta}^{\sigma'_{par(\sigma')}} \\
&= \sum_{i\in\Z}x_{i+1}^{\sigma_1}\dots x_{i+par(\sigma)}^{\sigma_{par(\sigma)}}y_{i+1+\Delta}^{\sigma'_1}\dots y_{i+par(\sigma')+\Delta}^{\sigma'_{par(\sigma')}}.
\end{align*}

Furthermore, $f^{(4)}_{\sigma,\sigma',\Delta}=f^{(4)}_{(\sigma'_{par(\sigma')},\dots,\sigma'_1),(\sigma_{par(\sigma)},\dots,\sigma_1),\Delta+par(\sigma')-par(\sigma)}$ for all $\sigma,\sigma',\Delta$.
Additionally, if $ord(\sigma)=0$ or $ord(\sigma')=0$, then $f^{(4)}_{\sigma,\sigma',\Delta}=f^{(4)}_{\sigma,\sigma',0}$ for all $\Delta\in\Z$.

\end{prop}

\subsection{The frieze group $F_5$} \label{sec_F5}

Remark that $F_5$ contains $F_2$ as a subgroup with the same induced action, and therefore the elements of $(P_k[X,Y])^{F_5}$ can be expressed as a formal infinite linear combination of $f^{(2)}_{\sigma,\sigma',\Delta}$.

Recall that $F_5\cong D_{\infty}\cong F_3$.
Equivalently to the beginning of \Cref{sec_F3}, we obtain that any element of $F_5$ can be written in the form $vg^z$ or $g^z$ for some $z\in\Z$.
As $g^2=t'$ and $vg=r'$ it follows that any element of $F_5$ can be expressed in the form ${t'}^z$, $v{t'}^z$, $g{t'}^z$ or $r'{t'}^z$ for some $z\in\Z$.

Similarly to \Cref{calc_v_onevar}, we calculate (with $m=par(\sigma)$ and $m'=par(\sigma')$) the action of $v$:

\begin{equation*}
\begin{split}
& v.x_{i+1}^{\sigma_1}x_{i+2}^{\sigma_2}\dots x_{i+m}^{\sigma_{m}}y_{i+1+\Delta}^{\sigma'_1}y_{i+2+\Delta}^{\sigma'_2}\dots y_{i+m'+\Delta}^{\sigma'_{m'}} \\
& = x_{j+1}^{\sigma_{m}}x_{j+2}^{\sigma_{m-1}}\dots x_{j+m}^{\sigma_1}y_{j+1+\Delta'}^{\sigma'_{m'}}y_{j+2+\Delta'}^{\sigma'_{m'-1}}\dots y_{j+m'+\Delta'}^{\sigma'_1}
\end{split}
\end{equation*}
where $j=-i-m-1$ and $\Delta'=m-m'-\Delta$.
A monomial is then preserved by $v{t'}^z$ if and only if $(\sigma_1,\dots,\sigma_{par(\sigma)})=(\sigma_{par(\sigma)},\dots,\sigma_1)$,  $(\sigma'_1,\dots,\sigma'_{par(\sigma')})=(\sigma'_{par(\sigma')},\dots,\sigma'_1)$ and $\Delta=0$ (remark that in that case $m$ is pair and therefore $-i-m-1$ and $i+1$ have the same parity).

With the same $m$ and $m'$, for $r'$ we obtain
\begin{equation*}
\begin{split}
& r'.x_{i+1}^{\sigma_1}x_{i+2}^{\sigma_2}\dots x_{i+m}^{\sigma_{m}}y_{i+1+\Delta}^{\sigma'_1}y_{i+2+\Delta}^{\sigma'_2}\dots y_{i+m'+\Delta}^{\sigma'_{m'}} \\
& = y_{j+1}^{\sigma_{m}}y_{j+2}^{\sigma_{m-1}}\dots y_{j+m}^{\sigma_1}x_{j+1+\Delta'}^{\sigma'_{m'}}x_{j+2+\Delta'}^{\sigma'_{m'-1}}\dots x_{j+m'+\Delta'}^{\sigma'_1}\\
& = x_{j'+1}^{\sigma'_{m'}}x_{j'+2}^{\sigma'_{m'-1}}\dots x_{j'+m'}^{\sigma'_1}y_{j'+1-\Delta'}^{\sigma_{m}}y_{j'+2-\Delta'}^{\sigma_{m-1}}\dots y_{j'+m'-\Delta'}^{\sigma_1}
\end{split}
\end{equation*}
where $j=-i-m-2$, $\Delta'=m-m'-\Delta$ and $j'=j-\Delta'=-i-2m-2+m'+\Delta$.
A monomial is then preserved by $r'{t'}^z$ if and only if $(\sigma_1,\dots,\sigma_{par(\sigma)})=(\sigma'_{par(\sigma')},\dots,\sigma'_1)$ and $j'$ is the same parity as $i$ (notice that in this case $-\Delta'=-m+m'+\Delta=\Delta$).
We have $j'=-i-m-2-\Delta'=-i-m-2+\Delta$, which is of the same parity as $i$ if and only if $m+\Delta$ is even.

\begin{prop}
For $k\geq1$, each element of $(P_k[X,Y])^{F_5}$ can be written in a unique way as a formal infinite linear combination on the set (see \Cref{composition})
\begin{align*}
&\left\{f^{(5)}_{\sigma,\sigma',\Delta}|\sigma,\sigma'\mbox{ are non-negative compositions with }ord(\sigma)+ord(\sigma')=k,\Delta\in\Z\right\}\\
&\bigcup\left\{f^{(5)'}_{\sigma,\sigma',\Delta}|\sigma,\sigma'\mbox{ are non-negative compositions with }ord(\sigma)+ord(\sigma')=k,\Delta\in\Z\right\}
\end{align*}
where for $(\sigma_1,\dots,\sigma_{par(\sigma)})\neq(\sigma_{par(\sigma)},\dots,\sigma_1)$ or  $(\sigma'_1,\dots,\sigma'_{par(\sigma')})\neq(\sigma'_{par(\sigma')},\dots,\sigma'_1)$ or $\Delta\neq0$, and also $(\sigma_1,\dots,\sigma_{par(\sigma)})\neq(\sigma'_{par(\sigma')},\dots,\sigma'_1)$ or $par(\sigma)+\Delta$ is odd:

\begin{equation*}
\begin{split}
f^{(5)}_{\sigma,\sigma',\Delta} & = \sum_{i\in\Z}{t'}^i.x_{1}^{\sigma_1}\dots x_{m}^{\sigma_{m}}y_{1+\Delta}^{\sigma'_1}\dots y_{m'+\Delta}^{\sigma'_{m'}}+v.\sum_{i\in\Z}{t'}^i.x_{1}^{\sigma_1}\dots x_{m}^{\sigma_{m}}y_{1+\Delta}^{\sigma'_1}\dots y_{m'+\Delta}^{\sigma'_{m'}} \\
& + g.\sum_{i\in\Z}{t'}^i.x_{1}^{\sigma_1}\dots x_{m}^{\sigma_{m}}y_{1+\Delta}^{\sigma'_1}\dots y_{m'+\Delta}^{\sigma'_{m'}}+r'.\sum_{i\in\Z}{t'}^i.x_{1}^{\sigma_1}\dots x_{m}^{\sigma_{m}}y_{1+\Delta}^{\sigma'_1}\dots y_{m'+\Delta}^{\sigma'_{m'}}
\end{split}
\end{equation*}
with $m=par(\sigma)$ and $m'=par(\sigma')$, and

\begin{equation*}
\begin{split}
f^{(5)'}_{\sigma,\sigma',\Delta} & = \sum_{i\in\Z}{t'}^i.x_{0}^{\sigma_1}\dots x_{m-1}^{\sigma_{m}}y_{\Delta}^{\sigma'_1}\dots y_{m'-1+\Delta}^{\sigma'_{m'}}+v.\sum_{i\in\Z}{t'}^i.x_{0}^{\sigma_1}\dots x_{m-1}^{\sigma_{m}}y_{\Delta}^{\sigma'_1}\dots y_{m'-1+\Delta}^{\sigma'_{m'}} \\
& + g.\sum_{i\in\Z}{t'}^i.x_{0}^{\sigma_1}\dots x_{m-1}^{\sigma_{m}}y_{\Delta}^{\sigma'_1}\dots y_{m'-1+\Delta}^{\sigma'_{m'}}+r'.\sum_{i\in\Z}{t'}^i.x_{0}^{\sigma_1}\dots x_{m-1}^{\sigma_{m}}y_{\Delta}^{\sigma'_1}\dots y_{m'-1+\Delta}^{\sigma'_{m'}}
\end{split}
\end{equation*}
which gives
\begin{equation*}
\begin{split}
f^{(5)}_{\sigma,\sigma',\Delta} & = \sum_{i\in\Z}x_{2i+1}^{\sigma_1}x_{2i+2}^{\sigma_2}\dots x_{2i+m}^{\sigma_{m}}y_{2i+1+\Delta}^{\sigma'_1}y_{2i+2+\Delta}^{\sigma'_2}\dots y_{2i+m'+\Delta}^{\sigma'_{m'}}\\
& + \sum_{i\in\Z}x_{2i-m}^{\sigma_{m}}x_{2i-m+1}^{\sigma_{m-1}}\dots x_{2i-1}^{\sigma_1}y_{2i-m-\Delta'}^{\sigma'_{m'}}y_{2i-m+1-\Delta'}^{\sigma'_{m'-1}}\dots y_{2i-1-\Delta}^{\sigma'_1}\\
& + \sum_{i\in\Z}x_{2i+2+\Delta}^{\sigma'_1}x_{2i+3+\Delta}^{\sigma'_2}\dots x_{2i+m'+1+\Delta}^{\sigma'_{m'}}y_{2i+2}^{\sigma_1}y_{2i+3}^{\sigma_2}\dots y_{2i+m+1}^{\sigma_{m}} \\
& + \sum_{i\in\Z}x_{2i-m'+\Delta+1}^{\sigma'_{m'}}x_{2i-m'+\Delta+2}^{\sigma'_{m'-1}}\dots x_{2i+\Delta}^{\sigma'_1}y_{2i-m'+\Delta+1+\Delta'}^{\sigma_{m}}\dots y_{2i+\Delta+\Delta'}^{\sigma_1}
\end{split}
\end{equation*}
with $\Delta'=\Delta+m'-m$; and similarly for $f^{(5)'}_{\sigma,\sigma',\Delta}$, but with opposite parity. And in all other cases $f^{(5)}_{\sigma,\sigma',\Delta}=f^{(2)}_{\sigma,\sigma',\Delta}$ and $f^{(5)'}_{\sigma,\sigma',\Delta}=f^{(2)'}_{\sigma,\sigma',\Delta}$ (see~\Cref{prop_f2}).

Furthermore, for any $\sigma,\sigma',\Delta$ we obtain the following equalities: If $par(\sigma)$ is odd we get $f^{(5)}_{\sigma,\sigma',\Delta}=f^{(5)}_{(\sigma_{m},\dots,\sigma_1),(\sigma'_{m'},\dots,\sigma'_1),m-m'-\Delta}$ (with $m=par(\sigma)$ and $m'=par(\sigma')$) and $f^{(5)'}_{\sigma,\sigma',\Delta}=f^{(5)'}_{(\sigma_{m},\dots,\sigma_1),(\sigma'_{m'},\dots,\sigma'_1),m-m'-\Delta}$.
If $par(\sigma)$ is even, we instead have $f^{(5)}_{\sigma,\sigma',\Delta}=f^{(5)'}_{(\sigma_{m},\dots,\sigma_1),(\sigma'_{m'},\dots,\sigma'_1),m-m'-\Delta}$.
Additionally, if $par(\sigma')+\Delta$ is even, we have $f^{(5)}_{\sigma,\sigma',\Delta}=f^{(5)}_{(\sigma'_{m'},\dots,\sigma'_1),(\sigma_{m},\dots,\sigma_1),\Delta+m'-m}$ and $f^{(5)'}_{\sigma,\sigma',\Delta}=f^{(5)'}_{(\sigma'_{m'},\dots,\sigma'_1),(\sigma_{m},\dots,\sigma_1),\Delta+m'-m}$.
Whereas if $par(\sigma')+\Delta$ is odd, $f^{(5)}_{\sigma,\sigma',\Delta}=f^{(5)'}_{(\sigma'_{m'},\dots,\sigma'_1),(\sigma_{m},\dots,\sigma_1),\Delta+m'-m}$.

Moreover, as in \Cref{prop_f2}, for $\sigma=\sigma'$ and $\Delta=0$ we get $f^{(5)}_{\sigma,\sigma,0}=f^{(5)'}_{\sigma,\sigma,0}$.
Finally, if $ord(\sigma)=0$ or $ord(\sigma')=0$, then $f^{(5)}_{\sigma,\sigma',\Delta}=f^{(5)}_{\sigma,\sigma',0}$ and $f^{(5)'}_{\sigma,\sigma',\Delta}=f^{(5)'}_{\sigma,\sigma',0}$ for all $\Delta\in\Z$.
\end{prop}

\subsection{The frieze group $F_6$} \label{sec_F6}

As this group contains the translation, we obtain the results of \Cref{t_twovar_eq}.
We calculate (with $m=par(\sigma)$ and $m'=par(\sigma')$):

\begin{align*}
&h.x_{i+1}^{\sigma_1}x_{i+2}^{\sigma_2}\dots x_{i+m}^{\sigma_{m}}y_{i+1+\Delta}^{\sigma'_1}y_{i+2+\Delta}^{\sigma'_2}\dots y_{i+m'+\Delta}^{\sigma'_{m'}}\\
&=y_{i+1}^{\sigma_1}y_{i+2}^{\sigma_2}\dots y_{i+m}^{\sigma_{m}}x_{i+1+\Delta}^{\sigma'_1}x_{i+2+\Delta}^{\sigma'_2}\dots x_{i+m'+\Delta}^{\sigma'_{m'}}\\
&=x_{j+1}^{\sigma'_1}x_{j+2}^{\sigma'_2}\dots x_{j+m'}^{\sigma'_{m'}}y_{j+1-\Delta}^{\sigma_1}y_{j+2-\Delta}^{\sigma_2}\dots y_{j+m-\Delta}^{\sigma_{m}}
\end{align*}
for $j=i+\Delta$.
A monomial is then fixed by $ht^z$ if and only if $\sigma=\sigma'$ and $\Delta=0$.

\begin{prop} \label{prop_f6}
For $k\geq1$, each element of $(P_k[X,Y])^{F_6}$ can be written in a unique way as a formal infinite linear combination on the set (see \Cref{composition})
$$\left\{f^{(6)}_{\sigma,\sigma',\Delta}|\sigma,\sigma'\mbox{ are non-negative compositions with }ord(\sigma)+ord(\sigma')=k,\Delta\in\Z\right\}$$
where for $\sigma\neq\sigma'$ or $\Delta\neq0$:
\begin{align*}
f^{(6)}_{\sigma,\sigma',\Delta} & = \sum_{i\in\Z}t^i.x_{1}^{\sigma_1}\dots x_{m}^{\sigma_{m}}y_{1+\Delta}^{\sigma'_1}\dots y_{m'+\Delta}^{\sigma'_{m'}}+h.\sum_{i\in\Z}t^i.x_{1}^{\sigma_1}\dots x_{m}^{\sigma_{m}}y_{1+\Delta}^{\sigma'_1}\dots y_{m'+\Delta}^{\sigma'_{m'}} \\
& = \sum_{i\in\Z}x_{i+1}^{\sigma_1}x_{i+2}^{\sigma_2}\dots x_{i+m}^{\sigma_{m}}y_{i+1+\Delta}^{\sigma'_1}y_{i+2+\Delta}^{\sigma'_2}\dots y_{i+m'+\Delta}^{\sigma'_{m'}}\\
& + \sum_{i\in\Z}x_{i+1}^{\sigma'_1}x_{i+2}^{\sigma'_2}\dots x_{i+m'}^{\sigma'_{m'}}y_{i+1-\Delta}^{\sigma_1}y_{i+2-\Delta}^{\sigma_2}\dots y_{i+m-\Delta}^{\sigma_{m}}
\end{align*}
where  $m=par(\sigma)$ and $m'=par(\sigma')$, and for $\sigma=\sigma'$ and $\Delta=0$:
\begin{align*}
f^{(6)}_{\sigma,\sigma,0}& = \sum_{i\in\Z}t^i.x_{1}^{\sigma_1}x_{2}^{\sigma_2}\dots x_{m}^{\sigma_{m}}y_{1}^{\sigma_1}y_{2}^{\sigma_2}\dots y_{m}^{\sigma_{m}} \\
& = \sum_{i\in\Z}x_{i+1}^{\sigma_1}x_{i+2}^{\sigma_2}\dots x_{i+m}^{\sigma_{m}}y_{i+1}^{\sigma_1}y_{i+2}^{\sigma_2}\dots y_{i+m}^{\sigma_{m}}.
\end{align*}

Furthermore, $f^{(6)}_{\sigma,\sigma',\Delta}=f^{(6)}_{\sigma',\sigma,-\Delta}$ for all $\sigma,\sigma',\Delta$.
Additionally, if $ord(\sigma)=0$ or $ord(\sigma')=0$, then $f^{(6)}_{\sigma,\sigma',\Delta}=f^{(6)}_{\sigma,\sigma',0}$ for all $\Delta\in\Z$.
\end{prop}

\subsection{The frieze group $F_7$} \label{sec_F7}

Remark that the actions on monomials of each generating element, as well as of the element $r$, have been described in previous sections.
Indeed, the actions of $t$ and $r$ are described in \Cref{sec_F4}, the action of $v$ in \Cref{sec_F5} and of $h$ in \Cref{sec_F6}.
Moreover, in those sections we also describe the conditions in which a monomial is fixed by some element of the group.
Therefore we can directly give the result:

\begin{prop}
For $k\geq1$, each element of $(P_k[X,Y])^{F_6}$ can be written in a unique way as a formal infinite linear combination on the set (see \Cref{composition})
$$\left\{f^{(7)}_{\sigma,\sigma',\Delta}|\sigma,\sigma'\mbox{ are non-negative compositions with }ord(\sigma)+ord(\sigma')=k,\Delta\in\Z\right\}$$
where for $(\sigma_1,\dots,\sigma_{par(\sigma)})\neq(\sigma'_{par(\sigma')},\dots,\sigma'_1)$, and also $(\sigma_1,\dots,\sigma_{par(\sigma)})\neq(\sigma_{par(\sigma)},\dots,\sigma_1)$ or $(\sigma'_1,\dots,\sigma'_{par(\sigma')})\neq(\sigma'_{par(\sigma')},\dots,\sigma'_1)$ or $\Delta\neq0$, and also $\sigma\neq\sigma'$ or $\Delta\neq0$:
\begin{equation*}
\begin{split}
f^{(7)}_{\sigma,\sigma',\Delta} & = \sum_{i\in\Z}t^i.x_{1}^{\sigma_1}\dots x_{m}^{\sigma_{m}}y_{1+\Delta}^{\sigma'_1}\dots y_{m'+\Delta}^{\sigma'_{m'}}+v.\sum_{i\in\Z}t^i.x_{1}^{\sigma_1}\dots x_{m}^{\sigma_{m}}y_{1+\Delta}^{\sigma'_1}\dots y_{m'+\Delta}^{\sigma'_{m'}} \\
& + h.\sum_{i\in\Z}t^i.x_{1}^{\sigma_1}\dots x_{m}^{\sigma_{m}}y_{1+\Delta}^{\sigma'_1}\dots y_{m'+\Delta}^{\sigma'_{m'}}+r.\sum_{i\in\Z}t^i.x_{1}^{\sigma_1}\dots x_{m}^{\sigma_{m}}y_{1+\Delta}^{\sigma'_1}\dots y_{m'+\Delta}^{\sigma'_{m'}}
\end{split}
\end{equation*}
with $m=par(\sigma)$ and $m'=par(\sigma')$, which gives
\begin{equation*}
\begin{split}
f^{(7)}_{\sigma,\sigma',\Delta} & = \sum_{i\in\Z}x_{i+1}^{\sigma_1}x_{i+2}^{\sigma_2}\dots x_{i+m}^{\sigma_{m}}y_{i+1+\Delta}^{\sigma'_1}y_{i+2+\Delta}^{\sigma'_2}\dots y_{i+m'+\Delta}^{\sigma'_{m'}}\\
& + \sum_{i\in\Z}x_{i+1}^{\sigma_{m}}x_{i+2}^{\sigma_{m-1}}\dots x_{i+m}^{\sigma_1}y_{i+1-\Delta'}^{\sigma'_{m'}}y_{i+2-\Delta'}^{\sigma'_{m'-1}}\dots y_{i+m'-\Delta'}^{\sigma'_1}\\
& + \sum_{i\in\Z}x_{i+1}^{\sigma'_1}x_{i+2}^{\sigma'_2}\dots x_{i+m'}^{\sigma'_{m'}}y_{i+1-\Delta}^{\sigma_1}y_{i+2-\Delta}^{\sigma_2}\dots y_{i+m-\Delta}^{\sigma_{m}} \\
& + \sum_{i\in\Z}x_{i+1}^{\sigma'_{m'}}x_{i+2}^{\sigma'_{m'-1}}\dots x_{i+m'}^{\sigma'_1}y_{i+1+\Delta'}^{\sigma_{m}}y_{i+2+\Delta'}^{\sigma_{m-1}}\dots y_{i+m+\Delta'}^{\sigma_1}
\end{split}
\end{equation*}
with $\Delta'=\Delta+par(\sigma')-par(\sigma)$.

And if $(\sigma_1,\dots,\sigma_{par(\sigma)})=(\sigma_{par(\sigma)},\dots,\sigma_1)$ or $(\sigma'_1,\dots,\sigma'_{par(\sigma')})=(\sigma'_{par(\sigma')},\dots,\sigma'_1)$ or $\Delta=0$ or $\sigma\neq\sigma'$, but $(\sigma_1,\dots,\sigma_{par(\sigma)})\neq(\sigma'_{par(\sigma')},\dots,\sigma'_1)$:
\begin{equation*}
\begin{split}
f^{(7)}_{\sigma,\sigma',\Delta} & = \sum_{i\in\Z}t^i.x_{1}^{\sigma_1}\dots x_{m}^{\sigma_{m}}y_{1+\Delta}^{\sigma'_1}\dots y_{m'+\Delta}^{\sigma'_{m'}}+v.\sum_{i\in\Z}t^i.x_{1}^{\sigma_1}\dots x_{m}^{\sigma_{m}}y_{1+\Delta}^{\sigma'_1}\dots y_{m'+\Delta}^{\sigma'_{m'}} \\
& = \sum_{i\in\Z}x_{i+1}^{\sigma_1}x_{i+2}^{\sigma_2}\dots x_{i+m}^{\sigma_{m}}y_{i+1+\Delta}^{\sigma'_1}y_{i+2+\Delta}^{\sigma'_2}\dots y_{i+m'+\Delta}^{\sigma'_{m'}}\\
& + \sum_{i\in\Z}x_{i+1}^{\sigma_{m}}x_{i+2}^{\sigma_{m-1}}\dots x_{i+m}^{\sigma_1}y_{i+1-\Delta'}^{\sigma'_{m'}}y_{i+2-\Delta'}^{\sigma'_{m'-1}}\dots y_{i+m'-\Delta'}^{\sigma'_1}
\end{split}
\end{equation*}
with $\Delta'=\Delta+par(\sigma')-par(\sigma)$.

And in the last case, which is $(\sigma_1,\dots,\sigma_{par(\sigma)})=(\sigma'_{par(\sigma')},\dots,\sigma'_1)$, $f^{(7)}_{\sigma,\sigma',\Delta}=f^{(6)}_{\sigma,\sigma',\Delta}$ (see \Cref{prop_f6}).

Furthermore, for all $\sigma,\sigma',\Delta$ we have $f^{(7)}_{\sigma,\sigma',\Delta}=f^{(7)}_{(\sigma_{m},\dots,\sigma_1),(\sigma'_{m'},\dots,\sigma'_1),m-m'-\Delta}=f^{(7)}_{\sigma',\sigma,-\Delta}=f^{(7)}_{(\sigma'_{m'},\dots,\sigma'_1),(\sigma_{m},\dots,\sigma_1),\Delta+m'-m}$ (with $m=par(\sigma)$ and $m'=par(\sigma')$).
Additionally, if $ord(\sigma)=0$ or $ord(\sigma')=0$, then $f^{(7)}_{\sigma,\sigma',\Delta}=f^{(7)}_{\sigma,\sigma',0}$ for all $\Delta\in\Z$.
\end{prop}

\section{Structure of $P_k[X]$ and $P_k[X,Y]$ as representations of $F_1$ and $F_6$} \label{sec_Structure}

We first aim at describing the structure of $P_k[X]$ as a module over $F_1$.

For each non-negative composition $\sigma$ of order $k$ we will define the subspace $M^{(1)}_\sigma$ of $P_k[X]$ as all formal infinite linear combinations of elements from the set $\{x_{i+1}^{\sigma_1}x_{i+2}^{\sigma_2}\dots x_{i+par(\sigma)}^{\sigma_{par(\sigma)}} | i\in\Z\}$.
Then $M^{(1)}_\sigma$ is closed under the action of $F_1$ and $P_k[X]$ as a module over $F_1$ splits as the following direct product 
 
$$P_k[X]=\prod_{\sigma|ord(\sigma)=k}M^{(1)}_\sigma.$$

Remark that for $k\geq2$ this is an infinite product as we can have an arbitrary number of parts in $\sigma$.
It is easy to check that each $M^{(1)}_\sigma$ is isomorphic to $K^\Z$ - indeed, it suffices to send the element $x_{i+1}^{\sigma_1}x_{i+2}^{\sigma_2}\dots x_{i+par(\sigma)}^{\sigma_{par(\sigma)}}$ into the infinite sequence with value $1$ at position $i$ and $0$ everywhere else.
We thus obtain that $P_k[X]$ as an $F_1$-module is, for $k\geq2$, isomorphic to an infinite countable product of copies of $K^\Z$, and $P_1[X]$ is isomorphic to $K^\Z$.

We also remark that $K^\ZZ$ as a module over $\ZZ$ has one one-dimensional submodule, spanned by the infinite sequence with 1 at every position. There is also one countable dimensional submodule consisting of all sequences with finitely many non-zero entries. The question of determining all uncountable dimensional $\ZZ$-submodules of $K^\ZZ$ is still an open problem.

Consider in contrast the action of $F_6$ on $P_k[X,Y]$. Here we obtain
$$P_k[X,Y]=\prod_{\sigma,\sigma',\Delta|ord(\sigma)+ord(\sigma')=k}M^{(6)}_{\sigma,\sigma',\Delta}$$
where $M^{(6)}_{\sigma,\sigma',\Delta}$ (with $m=par(\sigma)$ and $m'=par(\sigma')$) are the formal infinite linear combinations on the sets $\{x_{i+1}^{\sigma_1}x_{i+2}^{\sigma_2}\dots x_{i+m}^{\sigma_{m}}y_{i+1+\Delta}^{\sigma'_1}y_{i+2+\Delta}^{\sigma'_2}\dots y_{i+m'+\Delta}^{\sigma'_{m'}}|i\in\Z\}$ and $\{x_{i+1}^{\sigma'_1}x_{i+2}^{\sigma'_2}\dots x_{i+m'}^{\sigma'_{m'}}y_{i+1-\Delta}^{\sigma_1}y_{i+2-\Delta}^{\sigma_2}\dots y_{i+m-\Delta}^{\sigma_{m}}|i\in\Z\}$.
As described in \Cref{sec_F6}, two types of structures arise depending on whether it is true that simultaneously $\sigma=\sigma'$ and $\Delta=0$.

If it is true, then $h$ acts trivially on the monomial $x_{i+1}^{\sigma_1}x_{i+2}^{\sigma_2}\dots x_{i+m}^{\sigma_{m}}y_{i+1}^{\sigma_1}y_{i+2}^{\sigma_2}\dots y_{i+m}^{\sigma_{m}}$ for every $i$, and the space $M^{(6)}_{\sigma,\sigma,0}$ is then isomorphic to $K^\Z$ for the same reason as $M^{(1)}_\sigma$.
Remark that this is only possible if $k$ is even.

If either $\sigma\neq\sigma'$ or $\Delta\neq0$, we have obtained that the elements described in the definition of $M^{(6)}_{\sigma,\sigma',\Delta}$ are pairwise distinct.
We then obtain an isomorphism into $(K\bigoplus K)^\Z$ by sending $x_{i+1}^{\sigma_1}\dots x_{i+m}^{\sigma_{m}}y_{i+1+\Delta}^{\sigma'_1}\dots y_{i+m'+\Delta}^{\sigma'_{m'}}$ into an element that has $(1,0)$ at $i$, and $x_{i+1}^{\sigma'_1}\dots x_{i+m'}^{\sigma'_{m'}}y_{i+1-\Delta}^{\sigma_1}\dots y_{i+m-\Delta}^{\sigma_{m}}$ into one that has $(0,1)$.

Therefore for even $k$, $P_k[X,Y]$ as a module over $F_6$ splits into an infinite countable number of copies of $(K\bigoplus K)^\Z$ and an infinite countable number of copies of $K^\Z$, whereas for odd $k\geq3$, we have an infinite countable number of copies of $(K \bigoplus K)^\Z$ and nothing else.
And for $k=1$, it is isomorphic to $(K \bigoplus K)^\Z$.

\end{document}